\documentclass[a4paper,10pt,leqno]{amsart}
\usepackage[english]{babel}
\usepackage[cp1251]{inputenc}
\usepackage[T2A]{fontenc}
\usepackage{amsmath}
\usepackage{amssymb}

\newcommand{\ds}{\displaystyle}

\newcommand{\bp}{\begin{proof}[Proof:]}
\newcommand{\ep}{\end{proof}}
\newcommand\md{\!\!\!\!\!\!\pmod }

\newcommand{\R}{\ensuremath{\mathbb{R}}}

\newcommand{\N}{\ensuremath{\mathbb{N}}}

\newtheorem{thm}{Theorem}
\newtheorem{lem}{Lemma}
\newtheorem{remm}{Remmark}

\begin{document}

\title[On the distribution of $\alpha p^2+\beta$ modulo one]
{On the distribution of $\alpha p^2+\beta$ modulo one for primes $p$ such that $p+2$
has no more two prime divisors}

\author{T.L.Todorova}

\date{}

\begin{abstract}
A classical problem in analytic number theory is to study the
distribution of fractional part $\alpha p^k+\beta,\,k\ge 1$ modulo 1, where $\alpha$ is irrational
and $p$ runs over the set of primes. For $k=2$ we consider the subsequence
generated by the primes $p$ such that $p+2$ is an almost-prime
(the existence of infinitely many such $p$ is another topical
result in prime number theory) and prove that its distribution has
a similar property.

\smallskip


\end{abstract}

\maketitle

\section{Introduction and statements of the result}
The famous prime twins conjecture states that there exist
infinitely many primes $p$ such that $p+2$ is a prime too. This
hypothesis is still unproved but in 1973 Chen
\cite{Chen} proved that there are infinitely many
primes $p$ for which $p+2=\mathrm{P}_2$. (As usual $\mathrm{P}_r$
denotes an integer with no more than $r$ prime factors, counted
according to multiplicity).

Let $\alpha$ be irrational real number and $||x||$ denote the distance from $x$
to the nearest integer. The distribution of fractional parts of the sequence $\alpha n^k$,
$\alpha\in \mathbb{R}\setminus\mathbb{Q}$
was first considered by Hardy, Littlewood \cite{HardyLitt} and Weyl \cite{Weyl}.
The problem of distribution of the fractional parts of $\alpha p^k$, where p denotes a prime,
first was considered by Vinogradov (see Chapter 11 of \cite{Vin4} for the case $k = 1$, \cite{Vin5} for $k\ge 2$),
who showed that for any real $\beta$ there are infinitely many primes $p$ such that
\begin{equation}\label{alphap}
 || \alpha p+\beta||<p^{-\theta},
\end{equation}
where $\theta=1/5-\varepsilon$, $\varepsilon>0$ is arbitrary small.
After that many authors improved the upper bound of the
exponent $\theta$. The best result
is given by Matom\"{a}ki \cite{Mato1} with $\theta=1/3-\varepsilon$.
Another interesting problem is the study of the distribution of
the fractional part of $\alpha p^k$ with $2\le k\le 12$,
such Baker and Harman \cite{HarmanBaker}, Wong \cite{Wong} etc.
For $2\le k\le 12$ the best result is due to Baker and Harman \cite{HarmanBaker}.

In \cite{TT} Todorova and Tolev considered the primes $p$
such that $||\alpha p+\beta||<p^{-\theta}$ and $p+2=P_r$ and prove existence
of such primes with $\theta =1/100$ and $r=4$. Later Matom\"{a}ki \cite{Mato1}
and San Ying Shi \cite{Shi} have shown that  this actually holds whit $p+2=P_2$ and
$\theta =1/1000$ and $\theta =1.5/100$ respectively.

In \cite{ShiWu} Shi and Wu proved existence of infinitely many primes $p$ such that $||\alpha p^2+\beta||<p^{-\theta}$
and $p+2=P_4$ with $0<\theta<2/375$. In 2021 Xue, Li and Zhang \cite{XLZ} improved
the result of Shi and Wu with $0<\theta<10/1561$.

In this paper we evaluate exponential sums over well-separable numbers and improve the results of
Shi, Wu and Xue, Li and Zhang.

We will say that $d$ is a well-separable number of level $D\ge 1$ if for any $H, S \ge 1$ with $HS = D$,
there are integers $h\le H$, $s\le S$ such that $d=hs$.

\begin{thm} \label{Mat2} Suppose $\alpha\in \mathbb{R}\setminus\mathbb{Q}$ satisfies
conditions
\begin{equation}\label{alfa}
  \bigg|\alpha- \frac{a}{q}\bigg|<\frac{1}{q^2},\quad a \in \mathbb{Z}, \,;\; q \in \mathbb{N }, \quad (a, q) = 1,\quad q \geq 1\,,
\end{equation}
$K$ and $D$ are definied by (\ref{us2}), $\lambda(d)$ are complex numbers defined for
$d\le D$,
\begin{equation}\label{ksi}
  \lambda(d)\ll \tau (d)\quad \hbox{and}\quad \lambda(d)\not =0\quad\hbox{if}\quad d\quad\hbox{is well-separable number of level}\quad D\,,
\end{equation}
$c(k)\ll 1$ are complex numbers, $0<|k|\le K$.
Then for any arbitrary small $\varepsilon>0$ and $b\in \mathbb{Z}$ for the sum
\begin{equation}\label{c}
   W(x)= \sum\limits_{d\le D}\lambda(d)
    \sum\limits_{1\le |k|\le K}c(k)
    \sum\limits_{\substack{ n\sim x \\ n\equiv b\,(d)}}e\big((\alpha n^2+\beta )k\big )\Lambda (n)
\end{equation}
we have
\begin{equation}\label{OcenkaW}
  W\ll x^{\varepsilon}\bigg(\frac{xK}{\Delta^{\frac{1}{32}}}+\frac{xK}{q^{\frac{1}{32}}}+
  x^{\frac{15}{16}}K^{\frac{31}{32}}q^{\frac{1}{32}} +\frac{x\Delta ^{\frac{1}{2}}K}{q^{\frac{1}{4}}}+
 x^{\frac{1}{2}}\Delta ^{\frac{1}{2}}K^{\frac{3}{4}}q^{\frac{1}{4}} \bigg)\,.
\end{equation}
\end{thm}
\begin{remm} It is obvious that the Theorem \ref{Mat2} is true if function $\lambda(d)$ is well-factorable.
\end{remm}
\begin{remm}
The result is non-trivial as soon as 
\begin{equation*}
  \max\{x^{\omega},\,\Delta K^2\}<q<\min\bigg\{\frac{x^2}{K^{31}},\,\frac{x^2}{\Delta ^2K^{3}}\bigg\}
\end{equation*}
for some arbitrarily small and fixed $\omega>0$. In this case the estimate
\begin{equation*}\label{ocWokonch}
   W\ll x^{\varepsilon}\bigg(\frac{xK}{\Delta^{\frac{1}{32}}}+
\frac{x\Delta^{\frac{1}{2}}K}{q^{\frac{1}{4}}}
+x^{\frac{15}{16}}K^{\frac{31}{32}}q^{\frac{1}{32}}+
x^{\frac{1}{2}}\Delta ^{\frac{1}{2}}K^{\frac{3}{4}}q^{\frac{1}{4}}
\bigg)\,.
\end{equation*}
is valid.
\end{remm}
\begin{lem}\label{Lemseq} Suppose $\alpha \in \mathbb{R}\backslash \mathbb{Q}$ satisfies
conditions (\ref{alfa}), sum $W(x)$ is defined by (\ref{c}),
$\lambda(d)$ are complex numbers defined for
$d\le D$ and satisfying (\ref{ksi}) and (\ref{us2}), $c(k)\ll 1$ are complex
numbers $0<|k|\le K$.
Then there exist
a sequence
$$
\{x_j\}_{j=1}^{\infty},\;\lim\limits_{j\to
\infty}x_j=\infty\,,
$$
such that
\begin{equation*}
    W(x_j)\ll x_j^{1-\omega}\,,\;\;\;\;j=1,\,2,\,3,\ldots\;
\end{equation*}
for any $\omega>0$.
\end{lem}

\begin{thm}\label{osnTh} Let $\alpha \in \mathbb{R}\backslash \mathbb{Q}$  satisfies
conditions (\ref{alfa}),
$\beta \in \mathbb{R}$
 and let
$$
0<~\theta< ~\dfrac{1}{1296}\,,
$$
Then there are
infinitely many primes $p$ satisfying $p+2=\mathrm{P}_2$ and such
that
\begin{equation}\label{alfabeta}
    ||\alpha p^2+\beta||<p^{-\theta }.
\end{equation}
\end{thm}

{\bf Acknowledgements:} This work was supported Sofia University
Scientific Fund, grant 80-10-99/2023.

\section{Notation}
Let $x$ be a sufficiently large real number and
$\theta$ and $\rho $ be real constants
satisfying
\begin{equation}\label{us1}
    0<\theta<\frac{1}{200},\qquad
    \rho>32\theta , \qquad \rho+4\theta<\frac{11}{54}, \qquad \rho +\theta >0
\, .
\end{equation}
We shall specify $\rho$ and $\theta$  latter. We put
\begin{equation}\label{us2}
    \delta =\delta (x)=x^{-\theta },\quad \,K=\delta ^{-1}\log ^2x,\quad
    \Delta=x^{\rho },\quad D=\frac{x^{1/2}}{\Delta K^4}.
\end{equation}
By $p$ and $q$ we always denote primes. As usual $\varphi (n),\,\mu (n),\,\Lambda
(n)$ denote respectively Euler's function, M\"{o}bius'
function and Mangoldt's function. We denote by $\tau _k(n)$ the
number of solutions of the equation $m_1m_2\ldots m_k=n$ in
natural numbers $m_1,\,\ldots,m_k$ and $\tau _2(n)=\tau (n)$. Let
$(m_1,\,\ldots \,,m_k)$ and $[m_1,\,\ldots \,,m_k]$ be the
greatest common divisor and the least common multiple of
$m_1,\,\ldots,m_k$ respectively. Instead of $m\equiv n\,\pmod {k}$
we write for simplicity $m\equiv n(k)$. As usual, $||y||$ denotes
 the distance from $y$ to the
nearest integer, $e(y)=e^{2\pi iy}$. For positive $A$ and $B$ we
write $A\asymp B$ instead of $A\ll B\ll A$ and  $k\sim K$ means $K/2 \le k < K$. The letter
$\varepsilon$ denotes an arbitrary small positive number, not the
same in all appearances. For example this convention allows us to
write $x^{\varepsilon }\log x\ll x^{\varepsilon }$.

\section{Some Lemmas}
\begin{lem} \label{Tau}
Let $ k,\, l,\,m,\, n \in \mathbb{N} $; $ X, \varepsilon \in \mathbb{R}$;
$X\ge 2$, $k\ge 2$ and $ \varepsilon >0$. Then
\begin{equation*}
\begin{array}{llll}
(\mathrm{i})  \qquad&\ds\sum\limits_{n\le X}  \big( \, \tau_k(n) \, \big)^l  \ll_{ { }_{k, l} }
                         X(\log X)^{k^l-1} \; ;
                         \qquad \qquad  \; \\[12pt]
(\mathrm{ii})  \qquad&\ds\tau_k(n)\ll_{ { }_{k, \varepsilon} } n^{\varepsilon} \; .
\end{array}
\end{equation*}
\end{lem}
\bp See \cite{VinBasic}, ch. 3.
\ep

\begin{lem} \label{Trsum1}
Let $ X\ge 1 $ and $a,\,d\in \mathbb{N}$. Then
\begin{equation*}
   \bigg|
\mathop{\sum}_{\substack{
                      n\le X \\
                      n\equiv a\md d }}
e(\alpha n)\bigg|
           \ll \min \left(\frac{X}{d},
                 \frac{1} { ||\alpha d||} \right)
\end{equation*}
\end{lem}
\bp See \cite{Kar}, ch.6, \S 2.
\ep

\begin{lem} \label{Trsum2}
Let $X ,\,Y \in \mathbb{R}\,;\; X,\,Y\ge 1$
and $\alpha $ satisfied conditions (\ref{alfa}).
Then
\begin{equation*}
  \sum_{n\le X} \, \min \Big( \frac{XY}{n},\, \frac{1} { ||\alpha n|| } \Big)
    \, \ll \, XY \,
      \Big( \, \frac{1} {q} +
                \frac{1} {Y} +
                \frac { q }{ XY } \, \Big)\log (2Xq)\; . \notag
\end{equation*}
\end{lem}
\bp See Lemma 2.2 from \cite{Va}, ch. 2,\S 2.1.
\ep

\begin{lem}\label{Mat} Let $\mu,\,\zeta\in\N$; $x,\,M,\,J\in \R^{+}$, $\alpha\in \mathbb{R}\setminus\mathbb{Q}$,
and $\alpha$ satisfy conditions (\ref{alfa}).
Then for every arbitrary small $\varepsilon>0$ the inequality
\begin{multline*}
  \sum\limits_{m\sim M}\tau_{\mu}(m)\sum\limits_{j\sim J}\tau_{\zeta}(j)
\min\bigg\{\frac{x}{m^2j},\,\frac{1}{||\alpha m^2j||}\bigg\}
\ll x^{\varepsilon}\bigg(MJ+\frac{x}{M^{\frac{3}{2}}}+\frac{x}{Mq^{\frac{1}{2}}}+\frac{x^{\frac{1}{2}}q^{\frac{1}{2}}}{M}\bigg)
\end{multline*}
is fulfilled.
\end{lem}
\bp See Lemma 8, \cite{Mato1}.
\ep

\begin{lem} If $d|P(z)$, $z<D^{1/2}$, $\lambda ^{\pm}$ are Rosser's weigts and either $\lambda^{+}(d)\not =0$ or
$\lambda^{-}(d)\not =0$ then $d$ is well-separable number.
\end{lem}
\bp See Lemma 12.16, \cite{FriIwa}
\ep

\begin{thm}\label{linSieve} Let $x>1$, $2\le z\le D^{1/2}$, $s=\dfrac{\log D}{\log z}$,
$\;\omega (d)$ be a multiplicative function, $0<\omega (p)<p$. If
\begin{equation*}
  \left\{
    \begin{array}{ll}
      \mathcal{A}_d=\dfrac{\omega (d)}{d}x+r(x,\,d) &\hbox{if}\quad \mu(d)\ne 0 \\
      \sum\limits_{z_1\le p<z_2}\dfrac{\omega (p)}{p}=\log\bigg(\dfrac{\log z_2}{\log z_1}\bigg)+
O\bigg(\dfrac{1}{\log z_1}\bigg), & z_2>z_1\ge 2
    \end{array}
  \right.
\end{equation*}
Then
\begin{multline*}
  xV(z)\bigg(f(s)+O\bigg(\frac{1}{(\log D)^{1/3}}\bigg)\bigg) \le S(\mathcal{A},\,z)
  \le xV(z)\bigg(F(s)+O\bigg(\frac{1}{(\log D)^{1/3}}\bigg)\bigg)\,,
\end{multline*}
$f(s),\, F(s)$ are
determined by the following differential-difference equaqions
\begin{equation*}
  \begin{split}
     F(s)&=\frac{2e^{\gamma}}{s},\,f(s)=0\quad \hbox{if}\quad 0<s\le 2\, \\
    (sF(s))'&=f(s-1),\quad (sf(s))'=F(s-1)\quad \hbox{if}\quad s>2\,,
   \end{split}
\end{equation*}
where $\gamma$  denote the Euler’s constant.
\end{thm}
\bp See Lemma 3 \cite{Iwa1}.
\ep
\section{Auxiliary results}

\begin{lem}\label{M3S2J} Let $\alpha\in \mathbb{R}\setminus\mathbb{Q}$ satisfied conditions (\ref{alfa}),
$M,\,S,\,J,\,x\in \mathbb{R^{+}}$, $x>M^3S^2J$ and $\mu,\,\sigma,\,\zeta\in [2,\,\infty)\cap\mathbb{N}$,
\begin{equation*}
  G=\sum\limits_{m\sim M}\tau_{\mu}(m)
\sum\limits_{s\sim S}\tau_{\sigma}(s)
\sum\limits_{j\sim J}\tau_{\zeta} (j)
  \min\bigg\{\frac{x}{m^3s^2j},\,\frac{1}{||\alpha m^3s^2j||}\bigg\}\,.
\end{equation*}
Then for any $\varepsilon>0$ the inequalities
\begin{equation}\label{Ocz1}
G \ll x^{\varepsilon} \bigg(MSJ+\dfrac{x}{M^{\frac{9}{4}}S}+\dfrac{x}{M^2S^{\frac{9}{8}}}+
\dfrac{x}{M^2 S q^{\frac{1}{8}}}+\dfrac{x^{\frac{7}{8}} q^{\frac{1}{8}}}{M^2S}\bigg)
\end{equation}
and
\begin{equation}\label{Ocz2}
     G \ll x^{\varepsilon} \bigg(MSJ+\dfrac{x}{M^{\frac{9}{4}}S^{\frac{3}{4}}}+\dfrac{x}{M^2S^{\frac{3}{4}}q^{\frac{1}{4}}}+
\dfrac{x^{\frac{3}{4}} q^{\frac{1}{4}}}{M^2S^{\frac{3}{4}}}\bigg)\,
\end{equation}
are fulfilled.
\end{lem}
\bp Our proof is similar to proof of Lemma 8, \cite{Mato1}.
Let
\begin{equation}\label{H}
  H = \dfrac{x}{M^3S^2J}\,.
\end{equation}
If $H \le 2$, then trivially from Lemma \ref{Tau} (ii) we get
\begin{equation}\label{Hsmall}
  G \ll x^{\varepsilon}MSJ\,.
\end{equation}
So we can assume that $H > 2$.
From Lemma \ref{Tau} (ii) it is obviously that
\begin{equation*}
  G\ll x^{\varepsilon}\sum\limits_{m\sim M}
\sum\limits_{s\sim S}
\sum\limits_{j\sim J}
  \min\bigg\{\frac{x}{m^3s^2j},\,\frac{1}{||\alpha m^3s^2j||}\bigg\}\,.
\end{equation*}
We apply the Fourier expansion to
function $\min\bigg\{\dfrac{x}{m^3s^2j},\,\dfrac{1}{||\alpha m^3s^2j||}\bigg\}$ and get
\begin{equation*}
  \min\bigg\{\frac{x}{m^3s^2j},\,\frac{1}{||\alpha m^3s^2j||}\bigg\}=
  \sum\limits_{0<|h|\le H^2}w(h)e(\alpha m^3s^2jh) +O(\log x)\,,
\end{equation*}
where
\begin{equation}\label{wh}
  w(h)\ll \min\bigg\{\log H,\,\frac{H}{|h|}\bigg\}\,.
\end{equation}
Then
\begin{equation}\label{GI}
  |G|\ll x^{\varepsilon}\sum\limits_{0<|h|\le H^2}|w(h)|\sum\limits_{s\sim S}
\sum\limits_{j\sim J}\bigg |\sum\limits_{m\sim M}e(\alpha m^3s^2jh)\bigg |+
  MSJ\log x\,.
\end{equation}
So if
\begin{equation*}
  G(H_0)=\sum\limits_{h\sim H_0}\sum\limits_{s\sim S}
\sum\limits_{j\sim J}
  \bigg|\sum\limits_{m\sim M}e(\alpha m^3s^2jh)\bigg|\,.
\end{equation*}
then using (\ref{wh}) we have
\begin{equation}\label{G}
  G\ll x^{\varepsilon}\bigg(MSJ+\max\limits_{1\le H_0\le H_1}G(H_0)+
\max\limits_{H_1< H_0\le H^2}\frac{H}{H_0}G(H_0)\bigg)\,.
\end{equation}
We shall evaluate the sum $G(H_0)$. Applying the Cauchy-Schwarz inequality we obtain
\begin{align*}
  G^2(H_0)  &\ll x^{\varepsilon}H_0JS \sum\limits_{h\sim H_0}\sum\limits_{s\sim S}\sum\limits_{j\sim J}
  \bigg|\sum\limits_{m\sim M}e(\alpha m^3s^2jh)\bigg|^2\,\\
&\ll x^{\varepsilon}H_0JS \sum\limits_{h\sim H_0}\sum\limits_{s\sim S}\sum\limits_{j\sim J}
  \sum\limits_{m_1\sim M}\sum\limits_{m_2\sim M}e\big(\alpha (m_1^3-m_2^3)s^2jh\big)\,.
\end{align*}
Substituting $m_1=m_2+t$, where $0\le |t|\le M$ we get
\begin{equation}\label{G2H0}
  G^2(H_0)  \ll x^{\varepsilon}\bigg(H_0^2J^2S^2M+H_0JSG_1(H_0)\bigg)\,,
\end{equation}
where
\begin{equation*}
  G_1(H_0)=\sum\limits_{h\sim H_0}\sum\limits_{s\sim S}\sum\limits_{j\sim J}
  \sum\limits_{0<|t|<M}\bigg|\sum\limits_{m_2\sim M}e\big(\alpha (3m_2^2t+3m_2t^2)s^2jh\big)\bigg|\,.
\end{equation*}
Applying again the Cauchy-Schwarz inequality we obtain
\begin{multline*}
  G_1^2(H_0) \ll H_0JSM
\sum\limits_{h\sim H_0}\sum\limits_{s\sim S}\sum\limits_{j\sim J}
  \sum\limits_{0<|t|<M}\sum\limits_{m_2\sim M}\\
\times\sum\limits_{m_3\sim M}e\big(\alpha (3(m_2^2-m_3^2)t+3(m_2-m_3)t^2)s^2jh\big)\,.
\end{multline*}
Substituting $m_2=m_3+\ell$, where $0\le |\ell|\le M$ we get
\begin{multline*}
  G_1^2(H_0) \ll  H_0^2J^2S^2M^3+\\
H_0JSM\sum\limits_{h\sim H_0}\sum\limits_{s\sim S}\sum\limits_{j\sim J}
  \sum\limits_{0<|t|<M}\sum\limits_{0<|\ell|<M}\bigg|\sum\limits_{m_3\sim M}e\big(6\alpha m_3t\ell s^2jh\big)\bigg|\,.
\end{multline*}
Let $u=6t\ell hj$. Then using Lemma \ref{Trsum1} we get
\begin{equation*}
  G_1^2(H_0) \ll H_0^2J^2S^2M^3 +
 H_0JSM\sum\limits_{u\le 24H_0JM^2}\tau_5(u)\sum\limits_{s\sim S}\bigg|\sum\limits_{m_3\sim M}e\big(\alpha m_3s^2u\big)\bigg|
\end{equation*}
\begin{equation}\label{Twoways}
 \; \ll H_0^2J^2S^2M^3 +
 H_0JSM\sum\limits_{u\ll H_0JM^2}\tau_5(u)\sum\limits_{s\sim S}\min\bigg\{\frac{H_0JS^2 M^3}{s^2u},\,\frac{1}{||\alpha s^2u||}\bigg\}\,.
\end{equation}
We will estimate the above sum in two ways.
Using Lemma \ref{Mat} we obtain
\begin{equation*}
   G_1^2(H_0)\ll x^{\varepsilon}\bigg(H_0^2J^2S^2M^3+H_0^2J^2S^{\frac{3}{2}}M^{4}+
\frac{H_0^2J^2S^2M^4}{q^{\frac{1}{2}}}+H_0^{\frac{3}{2}}J^{\frac{3}{2}}SM^{\frac{5}{2}}q^{\frac{1}{2}}\bigg)\,.
\end{equation*}
So from (\ref{G2H0})
\begin{multline}\label{GH0}
  G(H_0)  \ll x^{\varepsilon}\bigg(H_0JSM^{\frac{3}{4}}+H_0JS^{\frac{7}{8}}M+
\frac{H_0JSM}{q^{\frac{1}{8}}}+
H_0^{\frac{7}{8}}J^{\frac{7}{8}}S^{\frac{3}{4}}
M^{\frac{5}{8}}q^{\frac{1}{8}}\bigg)\,.
\end{multline}
Choosing $H_0=H$ from (\ref{H}), (\ref{G}) and (\ref{GH0}) we get (\ref{Ocz1}).

On the other hand we can write the inequality (\ref{Twoways}) as
\begin{equation*}
  G_1^2(H_0) \ll H_0^2J^2S^2M^3 +
 H_0JSM\sum\limits_{k\ll H_0JS ^2M^2}\min\bigg\{\frac{H_0JS^2 M^3}{k},\,\frac{1}{||\alpha k||}\bigg\}\,
\end{equation*}
and using Lemma \ref{Trsum2} and (\ref{G2H0}) we get
\begin{equation}\label{Gh02}
  G^2(H_0) \ll x^{\varepsilon} \bigg(H_0^2J^2S^3M^3 +
 \frac{H_0^2J^2S^3M^4}{q}+H_0JSMq\bigg)\,
\end{equation}
Now we choose $H_0=H$. Then from (\ref{H}), (\ref{G}) and (\ref{Gh02})
the inequality (\ref{Ocz2}) is received.
\ep

\section{Proof of Theorem \ref{Mat2}}

To prove Theorem \ref{Mat2} we shall use the Vaughan's identity and we evaluate the sum $W$ in cases:
\begin{equation*}
  \begin{split}
    &\hbox{when } x^{\frac{8}{27}}\Delta\le D\le \dfrac{x^{\frac{1}{2}}}{\Delta K^4}\label{GraniciD1},  \hbox{and ;} \\
    &\hbox{when } D\le x^{\frac{8}{27}}\Delta\,.
  \end{split}
\end{equation*}

\subsection{Evaluation at large values of $D$}\label{Vaughan}
$\\$
Let
\begin{equation}\label{OgrD1}
  x^{\frac{8}{27}}\Delta\le D\le \dfrac{x^{\frac{1}{2}}}{\Delta K^4}
\end{equation}
First we decompose the sum $W(x)$ into $O(\log^2 x)$ sums of
type
\begin{equation*}\label{ccc}
   W=W(x,\,D,\,K)= \sum\limits_{d\sim D}\lambda(d)
    \sum\limits_{k\sim K_0}c(k)
    \sum\limits_{\substack{ n\sim x \\ n\equiv b\,(d)}}e\big((\alpha n^2+\beta )k\big )\Lambda (n)\,,
\end{equation*}
where 
\begin{equation}\label{OgrK0}
  1<K_0\le K
\end{equation}
and $\lambda(d)$ is Roser weight. As a necessary 
condition for $\lambda(d)\ne 0$ is numbers $d$ are 
squarefree from here on we will consider squarefree numbers $d$.
By Vaughan's identity we can decompose the sum $W$ into $O(\log x)$
type I sums
\begin{equation*}
  W_1=\mathop{\sum}_{\substack{
                      d\sim D \\
                      (b,\,d)=1}}\lambda(d)
\sum\limits_{k\sim K_0}c(k)e(\beta k)
\mathop{\sum}_{\substack{
                      m\sim M  \\
                      \ell\sim L\\
                      m\ell\equiv b\md d }}a(m)e(\alpha (m\ell)^2k)
\end{equation*}
or
\begin{equation*}
  W'_1=\mathop{\sum}_{\substack{
                      d\sim D \\
                      (b,\,d)=1}}\lambda(d)
\sum\limits_{k\sim K_0}c(k)e(\beta k)
\mathop{\sum}_{\substack{
                      m\sim M  \\
                      \ell\sim L\\
                      m\ell\equiv b\md d }}\log(\ell)e(\alpha (m\ell)^2k)
\end{equation*}
with $M\le x^{1/3}$ and $O(\log x)$ type II sums
\begin{equation}\label{ExpressionW2}
  W_2=\mathop{\sum}_{\substack{
                      d\sim D \\
                      (b,\,d)=1}}\lambda(d)
\sum\limits_{k\sim K_0}c(k)e(\beta k)
\mathop{\sum}_{\substack{
                      m\sim M  \\
                      \ell\sim L\\
                      m\ell\equiv b\md d }}a(m)b(\ell)e(\alpha (m\ell)^2k)
\end{equation}
with $M\in[x^{1/3},\,x^{2/3}]$ and

\begin{equation}\label{MNx}
  ML\sim x,\quad a(m)\ll \tau (m)\log m,\quad b(\ell)\ll \tau (\ell)\log\ell
\end{equation}

\subsubsection{Evaluation of type II sums.}
$\\$
The proof uses the idea of the proof of Theorem 1, \cite{Mato1}. As $x^{1/3}\le M,\,L\le x^{2/3}$ and
$ML\sim x$ we will consider only the case $x^{1/2}\le M\le x^{2/3}$.
The evaluation in the case $x^{1/2}\le L\le x^{2/3}$ is the same.
Using that $d$ is well-separable numbers we write
$d=hs$, where $(h,s)=1$ as $d$ is squarefree.  So the sum $W_2$ is presented as
$O(\log^2 x)$ sums of the type
\begin{align*}
  W_2=&\mathop{\sum}_{\substack{
                      h\sim H \\
                      (h,\,b)=1 }}
\mathop{\sum}_{\substack{
                      s\sim S   \\
                      (s,\,bh)=1 }}\lambda(hs)
\sum\limits_{k\sim K_0}c(k)e(\beta k)
\mathop{\sum}_{\substack{
                      \ell\sim L}}
\mathop{\sum}_{\substack{
                      m\sim M  \\
                      m\ell\equiv b\md {hs} }}a(m)b(\ell)e(\alpha (m\ell)^2k)\,.
\end{align*}
Here
\begin{equation}\label{HSD}
  h\sim H, \quad s\sim S, \quad D\sim HS\,
\end{equation}
and $H$ we will choose later.
Applying the Cauchy-Schwarz inequality to $W_{2}$
and using and Lemma \ref{Tau}(i) we obtain that
\begin{multline}\label{V0}
  W_2^2\ll x^{\varepsilon}K_0HM\sum\limits_{k\sim K_0}
\mathop{\sum}_{\substack{
                      h\sim H \\
                      (h,\,b)=1 }}
\mathop{\sum}_{\substack{
                      s'_1\sim S \\
                      (s'_1,\,bh)=1 }}\lambda(hs'_1)
\mathop{\sum}_{\substack{
                      s'_2\sim S \\
                      (s'_2,\,bh)=1 }}\lambda(hs'_2)\\
\times\sum\limits_{\ell_1\sim L}b(\ell_1)\sum\limits_{\ell_2\sim L}b(\ell_2)
\mathop{\sum}_{\substack{
                      m\sim M \\
                      m\ell_1\equiv b(hs'_1)\\
                      m\ell_2\equiv b(hs'_2)}}e(\alpha m^2(\ell_1^2-\ell_2^2)k)\,.
\end{multline}
Let $(s'_2,\,s'_1)=r$, $s'_1=rs_1$, $s'_2=rs_2$, $r\sim R,\,R\le S$ and $s'_1,\,s'_2\sim S/R$. Then
\begin{multline}\label{V1}
  W_2^2\ll x^{\varepsilon}K_0HM\sum\limits_{k\sim K_0}
\mathop{\sum}_{\substack{
                      h\sim H \\
                      (h,\,b)=1 }}
\mathop{\sum}_{\substack{
                      r\sim R \\
                      (r,\,bh)=1 }}
\mathop{\sum}_{\substack{
                      s_1\sim S/R \\
                      (s_1,\,bh)=1 }}\lambda(hrs_1)
\mathop{\sum}_{\substack{
                      s_2\sim S/R \\
                      (s_2,\,bs_1h)=1 }}\lambda(hrs_2)\\
\times\sum\limits_{\ell_1\sim L}b(\ell_1)\sum\limits_{\ell_2\sim L}b(\ell_2)
\mathop{\sum}_{\substack{
                      m\sim M \\
                      m\ell_1\equiv b(hrs_1)\\
                      m\ell_2\equiv b(hrs_2)}}e(\alpha m^2(\ell_1^2-\ell_2^2)k)\\
                  =W_{21}+W_{22},
\end{multline}
where $W_{21}$ is this one part of above sum for which
\begin{align*}
  \ell_1=\ell_2,\quad &s_1\ne s_2\quad\hbox{or} \\
  \ell_1=\ell_2,\quad &s_1=s_2=1\quad r\sim S\quad \hbox{or} \\
  \ell_1\ne \ell_2,\quad & s_1=s_2=1\quad r\sim S\quad\hbox{or} \\
  \ell_1\ne \ell_2,\quad &s_1\ne s_2,\quad M<\frac{4HS^2}{R}\,,
\end{align*}
$W_{22}$ is the rest part of sum for $W_2^2$.
Using Lemma \ref{Tau} and (\ref{OgrK0}) we get
\begin{equation}\label{W21}
  W_{21}\ll x^{\varepsilon}\bigg(xMHK^2+\frac{x^2HK^2}{D}+\frac{xLD^2K^2}{H}\bigg)\,.
\end{equation}
It is clear that for sum $W_{22}$ we have $\ell_1\ne \ell_2,\,s_1\ne s_2,\,M>\frac{4HS^2}{R}$.
From
$$
m\ell_1\equiv a(hrs_1),\; m\ell_2\equiv a(hrs_2)\quad \hbox{follows that}\quad \ell_1\equiv \ell_2 (hr)\,.
$$
We apply again the Cauchy-Schwarz inequality and get
\begin{align}\label{wrazkaW222}
  W_{22}^2&\ll  \frac{x^{2+\varepsilon}D^2K_0^3}{R^2}
\sum\limits_{k\sim K_0}
\mathop{\sum}_{\substack{
                      h\sim H \\
                      (h,\,b)=1 }}
\mathop{\sum}_{\substack{
                      r\sim R\\
                      (r,\,bh)=1 }}
\mathop{\sum}_{\substack{
                      s_1\sim S/R \\
                      (s_1,\,bh)=1 }}
\mathop{\sum}_{\substack{
                      s_2\sim S/R \\
                      (s_2,\,bs_1h)=1 }}\\
&\quad\quad\times\sum\limits_{\ell_1\sim L}\sum_{\substack{\ell_2\sim L\\\ell_2\equiv \ell_1 (hr)}}
\mathop{\sum}_{\substack{
                      m_1\sim M \\
                      m_1\ell_1\equiv b(hrs_1)\\
                      m_1\ell_2\equiv b(hrs_2) }}
\mathop{\sum}_{\substack{
                      m_2\sim M  \\
                      m_2\ell_1\equiv b(hrs_1)\\
                      m_2\ell_2\equiv b(hrs_2) }}e(\alpha (m_1^2-m_2^2)(\ell_1^2-\ell_2^2)k)\notag\\
&\quad\quad =W_{221}+W_{222}\,,\notag
\end{align}
where $W_{221}$ is this part of above sum for which $m_1=m_2$ and $W_{222}$ is the rest -
that is the sum which $m_1\ne m_2$. It is not difficult to see that
\begin{equation}\label{W221}
  W_{221}\ll \frac{x^{3+\varepsilon}LD^2K^4}{H}\,.
\end{equation}
Let consider the sum $W_{222}$. As
$$
m_i\ell_1\equiv b\pmod {hrs_1}\quad \hbox{and}\quad m_i\ell_2\equiv b\pmod {hrs_2},\,i=1,2
$$
we get
\begin{equation*}
  m_1\equiv m_2\md {hrs_1s_2}\equiv f\md {hrs_1s_2},\quad\hbox{where} \quad f=f(h,\,r,\,s_1,\,s_2,\,\ell_1,\,\ell_2)\,
\end{equation*}
and $\ell_1\equiv \ell_2\;\,\md {hr}$. Let
\begin{equation*}
  m_1=m_2+hrs_1s_2t,\,0<|t|\le \frac{8MR}{HS^2}\quad\hbox{and}\quad \ell_1=\ell_2+hru,\,0<|u|\le \frac{2L}{HR}\,.
\end{equation*}
Then
\begin{equation*}
  m_1^2-m_2^2=2m_2hrs_1s_2t+h^2r^2s_1^2s_2^2t^2\quad\hbox{and}\quad \ell_1^2-\ell_2^2=hru(2\ell_2+hru)\,.
\end{equation*}
So using above equalities and Lemma \ref{Trsum1} we obtain
\begin{align*}
  W_{222}&\ll  \frac{x^{2+\varepsilon}D^2K_0^3}{R^2}
\sum\limits_{k\sim K_0}
\mathop{\sum}_{\substack{
                      h\sim H }}
\mathop{\sum}_{\substack{
                      r\sim R }}
\mathop{\sum}_{\substack{
                      s_1\sim S/R  }}
\mathop{\sum}_{\substack{
                      s_2\sim S/R  }}\\
&\quad\quad\times\sum\limits_{\ell\sim L}
\sum\limits_{0<|u|\le \frac{2L}{HR}}
\sum\limits_{0<|t|\le \frac{8MR}{HS^2}}
\min\bigg\{\frac{M}{hrs_1s_2},\frac{1}{||2\alpha h^3r^3s_1^2s_2^2tu\ell k||}\bigg\},
\end{align*}
where $\ell=2\ell_2+hru$. We put
$$
m=hr,\quad s=s_1s_2,\quad j=2tunk,\;j\ll \dfrac{xLK_0}{D^2}
$$
and it is clear that the sum $W_{222}$ can be represented as a finite number of sums of the type
\begin{align*}
W_{223}=  \frac{x^{2+\varepsilon}D^2 K_0^3}{R^2}
\sum\limits_{m\sim HR}\tau (m)
\sum\limits_{s\sim \frac{S^2}{R^2}}\tau (s)
\sum\limits_{j\ll \frac{xLK_0}{D^2}}\tau_5(j)
\min\bigg\{\frac{x^2K_0}{m^3s^2j},\frac{1}{||\alpha m^3s^2j||}\bigg\}\,,
\end{align*}
so
\begin{equation}\label{W222W3}
  W_{222}\ll W_{223}\,.
\end{equation}
Using the inequality (\ref{Ocz1}) from Lemma \ref{M3S2J}, (\ref{HSD}) and (\ref{OgrK0}) we get
\begin{equation}\label{W223}
W_{223}\ll  x^{\varepsilon}\bigg(\frac{x^{3}LD^2K^4}{H}+
\frac{x^{4}K^4}{H^{\frac{1}{4}}}+\frac{x^4H^{\frac{1}{4}}K^4}{D^{\frac{1}{4}}}
+\frac{x^4K^4}{q^{\frac{1}{8}}}+x^{\frac{15}{4}}K^{\frac{31}{8}} q^{\frac{1}{8}}\bigg).
\end{equation}
From (\ref{V1}), (\ref{wrazkaW222}) and (\ref{W222W3}) we have
\begin{equation*}
  W_2\ll x^{\varepsilon}\bigg(W_{21}^{\frac{1}{2}}+W_{221}^{\frac{1}{4}}+W_{223}^{\frac{1}{4}}\bigg)
\end{equation*}
and from above inequality, (\ref{W21}), (\ref{W221}) and (\ref{W223}) we receive
\begin{multline}\label{W2_pyrvonachalno}
W_{2}\ll  x^{\varepsilon}\bigg(x^{\frac{1}{2}}M^{\frac{1}{2}}H^{\frac{1}{2}}K+\frac{x^{\frac{1}{2}}L^{\frac{1}{2}}DK}{H^{\frac{1}{2}}}+
\frac{x^{\frac{3}{4}}L^{\frac{1}{4}}D^{\frac{1}{2}}K}{H^{\frac{1}{4}}}\\
+\frac{xK}{H^{\frac{1}{16}}}+\frac{xH^{\frac{1}{16}}K}{D^{\frac{1}{16}}}+
\frac{xK}{q^{\frac{1}{32}}}+x^{\frac{15}{16}}K^{\frac{31}{32}} q^{\frac{1}{32}}\bigg).
\end{multline}
According to $D,\,M$ and $L$ we have
\begin{equation}\label{W2chrezV}
  W_{2}\ll x^{\varepsilon}\bigg(V_{1}+V_{2}+V_{3}+V_{4}\bigg)\,,
\end{equation}
where $V_{1}$ is the sum with
\begin{equation}\label{V1usl}
\begin{split}
        &x^{\frac{1}{2}}\le M\le \frac{x}{D},  \quad x^{\frac{2}{5}}\ll D\ll \frac{x^{\frac{1}{2}}}{\Delta K^4} \\
        &D\le L\le x^{\frac{1}{2}} \,,
\end{split}
\end{equation}
$V_{2}$ is the sum with
\begin{equation}\label{V1us2}
  \begin{split}
        &\frac{x}{D}\le M < x^{\frac{1}{3}}D^{\frac{2}{3}},  \quad x^{\frac{2}{5}}\ll D\ll \frac{x^{\frac{1}{2}}}{\Delta K^4} \\
        &\frac{x^{\frac{2}{3}}}{D^{\frac{2}{3}}}< L\le D  \,,
\end{split}
\end{equation}
$V_{3}$ is the sum with
\begin{equation}\label{V1us3}
  \begin{split}
        &x^{\frac{1}{3}}D^{\frac{2}{3}}\le M\le x^{\frac{2}{3}},  \quad x^{\frac{2}{5}}\ll D\ll \frac{x^{\frac{1}{2}}}{\Delta K^4} \\
        &x^{\frac{1}{3}}\le L\le \frac{x^{\frac{2}{3}}}{D^{\frac{2}{3}}}
\end{split}
\end{equation}
and $V_{4}$ is the sum with
\begin{equation}\label{V1us4}
  \begin{split}
        &x^{\frac{1}{2}}\le M\le x^{\frac{2}{3}},  \quad x^{\frac{8}{27}}\Delta\ll D\ll x^{\frac{2}{5}} \\
        &x^{\frac{1}{3}}\le L\le x^{\frac{1}{2}}  \,.
\end{split}
\end{equation}
For sums $V_1$, $V_2$, $V_3$ and $V_4$ we choose consequently
\begin{equation*}
  H=\frac{D}{\Delta^{1/2}},\quad H=\frac{LD^{2/3}}{x^{1/3}},\quad H=\frac{x^{1/3}}{\Delta },
\quad H=\frac{L^{4/5}D^{9/5}}{x^{4/5}}\,.
\end{equation*}
and from (\ref{W2_pyrvonachalno}), (\ref{W2chrezV}), (\ref{V1usl}), (\ref{V1us2}), (\ref{V1us3}) and (\ref{V1us4}) we get
\begin{equation*}
  W_{2}\ll \left\{
             \begin{array}{ll}
               x^{\varepsilon}\bigg(\dfrac{xK}{\Delta^{\frac{1}{32}}}
+\dfrac{xK}{q^{\frac{1}{32}}}+x^{\frac{5}{6}}D^{\frac{1}{3}}\Delta^{\frac{1}{4}}K
  +x^{\frac{15}{16}}K^{\frac{31}{32}} q^{\frac{1}{32}}\bigg), & \hbox{if}\quad x^{\frac{2}{5}}\ll D\ll \dfrac{x^{\frac{1}{2}}}{\Delta K^4}\,,\\
              x^{\varepsilon}\bigg(\dfrac{x^{\frac{31}{30}}K}{D^{\frac{9}{80}}}+\dfrac{xK}{q^{\frac{1}{32}}}
  +x^{\frac{15}{16}}K^{\frac{31}{32}} q^{\frac{1}{32}}\bigg), & \hbox{if}\quad x^{\frac{8}{27}}\Delta\ll D\ll x^{\frac{2}{5}}\,.
             \end{array}
           \right.
\end{equation*}
So
\begin{equation}\label{ocaW2}
  W_{2}\ll x^{\varepsilon}\bigg(\dfrac{xK}{\Delta^{\frac{1}{32}}}
+\dfrac{xK}{q^{\frac{1}{32}}}+x^{\frac{5}{6}}D^{\frac{1}{3}}\Delta^{\frac{1}{4}}K
  +x^{\frac{15}{16}}K^{\frac{31}{32}} q^{\frac{1}{32}}\bigg)\,.
\end{equation}
This estimate is nontrivial for
\begin{equation}\label{Ogr_q1}
  q<\frac{x^2}{K^{31}}\quad\hbox{и}\quad q>x^{\omega}
\end{equation}
for some arbitrarily small and fixed $\omega>0$.
\subsubsection{Evaluation of type I sums.}
$\\$
In this case we have that $L>x^{\frac{2}{3}}$ and $M<x^{\frac{1}{3}}$.
Again we will use that $d$ is well-separated numbers. So we can write
$d=hs$ whith $h$ and $s$ satisfying conditions (\ref{HSD}) 
and we will choose $H$ later. So the sum $W_1$ is presented as
$O(\log^2 x)$ sums of the type
\begin{align*}
  W_1=&\mathop{\sum}_{\substack{
                      h\sim H \\
                      (h,\,b)=1 }}
\mathop{\sum}_{\substack{
                      s\sim S   \\
                      (s,\,b)=1 }}\lambda(hs)
\sum\limits_{k\sim K_0}c(k)e(\beta k)
\mathop{\sum}_{\substack{
                      \ell\sim L }}
\mathop{\sum}_{\substack{
                      m\ell\sim M  \\
                      m\ell\equiv b\md d }}a(m)e(\alpha (m\ell)^2k)\,.
\end{align*}
Working in the same way as in the evaluation of the sum $W_2$ see (\ref{V0}), we get
\begin{multline}\label{V00}
  W_1^2\ll x^{\varepsilon}K_0HM\sum\limits_{k\sim K_0}
\mathop{\sum}_{\substack{
                      h\sim H \\
                      (h,\,b)=1 }}
\mathop{\sum}_{\substack{
                      r\sim R \\
                      (r,\,bh)=1 }}
\mathop{\sum}_{\substack{
                      s_1\sim S/R \\
                      (s_1,\,bh)=1 }}\lambda(hrs_1)
\mathop{\sum}_{\substack{
                      s_2\sim S/R \\
                      (s_2,\,bs_1h)=1 }}\lambda(hrs_2)\\
\times\sum\limits_{\ell_1\sim L}b(\ell_1)\sum\limits_{\ell_2\sim L}b(\ell_2)
\mathop{\sum}_{\substack{
                      m\sim M \\
                      m\ell_1\equiv b(hrs_1)\\
                      m\ell_2\equiv b(hrs_2)}}e(\alpha m^2(\ell_1^2-\ell_2^2)k)\\=x^{\varepsilon}\big(W_{11}+W_{12}+W_{13}\big)\,,
\end{multline}
where $W_{12}$ is this one part of above sum for which $\ell_1=\ell_2$,
$W_{13}$ is this one part of above sum for which $\ell_1\ne \ell_2,\quad  s_1=s_2=1\quad r\sim S$,
and $W_{11}$ is the rest part of sum for $W_1^2$.
Using that $M<x^{\frac{1}{3}}$ and (\ref{OgrK0}) we get
\begin{equation}\label{ocW12}
  W_{12}\ll x^{\varepsilon}.xMHK^2\,.
\end{equation}
Fot the sum $W_{13}\label{W13}$ we get
\begin{equation*}
  W_{13}\ll x^{\varepsilon}K_0HM\sum\limits_{k\sim K_0}
\mathop{\sum}_{\substack{
                      d\sim D \\
                      (d,\,b)=1 }}
\mathop{\sum}_{\substack{
                      m\sim M }}
\bigg |\mathop{\sum}_{\substack{
                      \ell_i\sim L\\
                      m\ell_i\equiv b(d)\\
                      i=1,1}}
                           e(\alpha m^2(\ell_1^2-\ell_2^2)k\bigg |\,.
\end{equation*}
As $\ell_1\equiv \ell_2\;\,\md {d}$ we put $\ell_1=\ell_2+du,\,0<|u|\ll \frac{L}{D}$.
So
\begin{equation*}
  W_{13}\ll x^{\varepsilon}K_0HM\sum\limits_{k\sim K_0}
\mathop{\sum}_{\substack{
                      d\sim D \\
                      (b,\,d)=1 }}
\mathop{\sum}_{\substack{
                      m\sim M }}
\mathop{\sum}_{\substack{
                      u\ll \frac{L}{D} }}
\bigg |\mathop{\sum}_{\substack{
                      \ell_2\sim L\\
                      m\ell_2\equiv b(d)}}
                           e(2\alpha m^2\ell_2udk)\bigg |\,
\end{equation*}
and from Lemma \ref{Trsum1} we get
\begin{align}\label{W13W14}
  W_{13}&\ll x^{\varepsilon}K_0HM\sum\limits_{k\sim K_0}
\mathop{\sum}_{\substack{
                      d\sim D  }}
\mathop{\sum}_{\substack{
                      m\sim M }}
\mathop{\sum}_{\substack{
                      u\ll \frac{L}{D} }}
\min\bigg\{\frac{x^2K_0}{m^2d^2(2uk)},\,\frac{1}{||\alpha m^2d^2(2uk)||}
\bigg\}\notag\\
&\ll x^{\varepsilon}K_0HM\sum\limits_{z\sim MD}\tau (z)
\mathop{\sum}_{\substack{
                      t\ll \frac{LK_0}{D} }}\tau_3(t)
\min\bigg\{\frac{x^2K_0}{z^2t},\,\frac{1}{||\alpha z^2t||}
\bigg\}\,.
\end{align}
Using Lemma \ref{Mat}, (\ref{W13W14}), $ML\sim x$ and (\ref{OgrK0}) we obtain
\begin{equation}\label{ocW13}
   W_{13}^{\frac{1}{2}}\ll x^{\varepsilon}\bigg(
x^{\frac{1}{2}}M^{\frac{1}{2}}H^{\frac{1}{2}}K+\frac{x^{\frac{3}{4}}L^{\frac{1}{4}}H^{\frac{1}{2}}K}{D^{\frac{3}{4}}}
+\frac{xH^{\frac{1}{2}}K}{D^{\frac{1}{2}}q^{\frac{1}{4}}}+\frac{x^{\frac{1}{2}}H^{\frac{1}{2}}q^{\frac{1}{4}}K^{\frac{3}{4}}}{D^{\frac{1}{2}}}
\bigg)
\end{equation}
Similarly we estimate the sum $W_{11}$ and get
\begin{multline*}\label{W11}
  W_{11}\ll x^{\varepsilon}K_0HM\sum\limits_{k\sim K_0}
\mathop{\sum}_{\substack{
                      h\sim H \\
                      (h,\,b)=1 }}
\mathop{\sum}_{\substack{
                      r\sim R \\
                      (r,\,bh)=1 }}
\mathop{\sum}_{\substack{
                      s_i\sim S/R \\
                      (s_1s_2,\,bh)=1 }}
\\
\times
\mathop{\sum}_{\substack{
                      m\sim M }}
\mathop{\sum}_{\substack{
                      u\ll \frac{L}{HR} }}
\bigg |\mathop{\sum}_{\substack{
                      \ell_2\sim L\\
                      m\ell_2\equiv b(hrs_2)\\
                      m(\ell_2+uhr)\equiv b(s_1)}}e(2\alpha m^2\ell_2uhrk)\bigg |\,.
\end{multline*}
So from Lemma \ref{Trsum1}
\begin{multline*}
   W_{11}\ll x^{\varepsilon}K_0HM\sum\limits_{k\sim K_0}
\mathop{\sum}_{\substack{
                      h\sim H }}
\mathop{\sum}_{\substack{
                      r\sim R }}
\mathop{\sum}_{\substack{
                      s_i\sim S/R \\
                      i=1,2}}
\\
\sum\limits_{\substack{
                      m\sim M }}
\mathop{\sum}_{\substack{
                      u\ll \frac{L}{HR} }}
\min\bigg\{\frac{L}{hrs_1s_2},\,\frac{1}{||\alpha m^2h^2r^2.2s_1s_2uk||}
\bigg\}\,
\end{multline*}
and
\begin{equation*}
   W_{11}\ll x^{\varepsilon}K_0HM\sum\limits_{z\sim MHR}\tau_3(z)
\mathop{\sum}_{\substack{
                      t\ll \frac{S^2LK}{HR^3} }}
\min\bigg\{\frac{x^2K_0}{z^2t},\,\frac{1}{||\alpha z^2t||}
\bigg\}\,,
\end{equation*}
where $z=mhr$, $t=2s_1s_2uk$.
Applying Lemma \ref{Mat} and using that $ML\sim x$ and (\ref{OgrK0}) we receive
\begin{equation}\label{ocW11}
W_{11}^{\frac{1}{2}}\ll x^{\varepsilon}\bigg(
\frac{x^{\frac{1}{2}}M^{\frac{1}{2}}DK}{H^{\frac{1}{2}}}
+\frac{x^{\frac{3}{4}}L^{\frac{1}{4}}K}{H^{\frac{1}{4}}}+\frac{xK}{q^{\frac{1}{4}}}+x^{\frac{1}{2}}q^{\frac{1}{4}}K^{\frac{3}{4}}
\bigg)
\end{equation}
Choosing $H=\frac{D}{M^{\frac{1}{2}}}$ from (\ref{V00}), (\ref{ocW12}), (\ref{ocW13}), (\ref{ocW11}) and (\ref{OgrD1}) follows
\begin{equation}\label{ocenkaW1}
  W_{1}\ll x^{\varepsilon}\bigg( x^{\frac{3}{4}}D^{\frac{1}{2}}K  +
\frac{xK}{D^{\frac{1}{4}}}+\frac{xK}{q^{\frac{1}{4}}}
  +x^{\frac{1}{2}}K^{\frac{3}{4}} q^{\frac{1}{4}}\bigg)\,.
\end{equation}
From (\ref{ocaW2}), (\ref{ocenkaW1}) and (\ref{Ogr_q1})
follows that in case  $x^{\frac{8}{27}}\Delta\le D\le \dfrac{x^{\frac{1}{2}}}{\Delta K^4}$
the estimate
\begin{equation}\label{ocenkaW3}
  W\ll x^{\varepsilon}\bigg(\dfrac{xK}{\Delta^{\frac{1}{32}}}
+\dfrac{xK}{q^{\frac{1}{32}}}
  +x^{\frac{15}{16}}K^{\frac{31}{32}} q^{\frac{1}{32}}\bigg)\,
\end{equation}
is fulfilled.
\subsection{Evaluation at small values of $D$\label{heath}}
$\\$
First we evaluate the second type sums.

\subsubsection{Evaluation of type II sums.}
$\\$
Let
\begin{equation*}
  D<x^{8/27}\Delta ,\quad x^{1/2}\le M\le x^{2/3},\quad x^{1/3}\le L\le x^{1/2}\,.
\end{equation*}
Applying the Cauchy-Schwarz inequality to $W_{2}$ (see (\ref{ExpressionW2}))
and using Lemma \ref{Tau}(i) and (\ref{OgrK0}) we obtain
\begin{multline}\label{VV0}
  W_2^2\ll x^{\varepsilon}K_0DM\sum\limits_{k\sim K_0}
\mathop{\sum}_{\substack{
                      d\sim D \\
                      (d,\,b)=1 }}
\sum\limits_{\substack{\ell_i\sim L\\ i=1,2\\ \ell_1\equiv \ell_2(d)}}
                      b(\ell_1)b(\ell_2)
\mathop{\sum}_{\substack{
                      m\sim M \\
                      m\ell_i\equiv b(d)\\
                      i=1,2}}e(\alpha m^2(\ell_1^2-\ell_2^2)k)\\
=V_{2}+x^{1+\varepsilon}MDK^2\,,
\end{multline}
where
\begin{equation*}
  V_{2}=x^{\varepsilon}K_0DM\sum\limits_{k\sim K_0}
\mathop{\sum}_{\substack{
                      d\sim D \\
                      (d,\,b)=1 }}
\sum\limits_{\substack{\ell_i\sim L\\
                        \ell_1\ne\ell_2\\
                        \ell_1\equiv \ell_2(d)}}b(\ell_1)b(\ell_2)
\mathop{\sum}_{\substack{
                      m\sim M \\
                      m\ell_i\equiv b(d)\\
                      i=1,2}}e(\alpha m^2(\ell_1^2-\ell_2^2)k)
\end{equation*}
Applying again the Cauchy-Schwarz inequality (working as in the evaluation of the sum $W_{222}$) we substitute
\begin{equation*}
  m_1=m_2+td,\quad t\ll \frac{M}{D}\quad\hbox{and}\quad \ell_1=\ell_2+\omega d,\quad \omega\ll\frac{L}{D}
\end{equation*}
and consequentially obtain
\begin{multline}\label{W12pred}
  V_{2}^2\ll x^{2+\varepsilon}D^2K_0^3\sum\limits_{k\sim K_0}
\mathop{\sum}_{\substack{
                      d\sim D \\
                      (d,\,b)=1 }}
\sum\limits_{\substack{\ell_i\sim L\\
                        \ell_1\ne\ell_2\\
                        \ell_1\equiv \ell_2(d)}}
\mathop{\sum}_{\substack{
                      m_i\sim M \\
                      m_1\ell_i\equiv b(d)\\
                      m_2\ell_i\equiv b(d)\\
                      m_1\ne m_2}}e(\alpha (m_1^2-m^2)(\ell_1^2-\ell_2^2)k)\\
+x^{3+\varepsilon}LD K^4\\
=V_{22}+x^{3+\varepsilon}LD K^4\,
\end{multline}
with
\begin{equation*}
   V_{22}\ll x^{2+\varepsilon}D^2K_0^3
\sum\limits_{k\sim K_0}
\mathop{\sum}_{\substack{
                      d\sim D }}
\mathop{\sum}_{\substack{
                      \omega\ll\frac{L}{D}}}
\mathop{\sum}_{\substack{
                      \ell_2\sim L}}
\mathop{\sum}_{\substack{
                      t\ll \frac{M}{D} }}
\min\bigg\{\frac{M}{d},\,\frac{1}{||\alpha d^3t\omega (2\ell_2+\omega d)k||}\,.
\bigg\}
\end{equation*}
Putting $\ell=\ell_2+\omega d$ and $z=\omega \ell kt$ we get
\begin{equation}\label{W22pred}
   V_{22}\ll x^{2+\varepsilon}D^2K_0^3
\mathop{\sum}_{\substack{
                      d\sim D }}
\mathop{\sum}_{\substack{
                      z\ll\frac{xLK_0}{D^2}}}
\tau_4(z)
\min\bigg\{\frac{x^2K_0}{d^3z},\,\frac{1}{||\alpha d^3z||}\,.
\bigg\}
\end{equation}
So from (\ref{VV0}), (\ref{W12pred}) and (\ref{OgrK0})
we receive
\begin{equation}\label{W2crezW22}
  W_{2}\ll x^{\varepsilon}\bigg(V_{22}^{\frac{1}{4}}+x^{\frac{1}{2}}M^{\frac{1}{2}}D^{\frac{1}{2}}K+x^{\frac{3}{4}}L^{\frac{1}{4}}D^{\frac{1}{4}}K\bigg)
\ll x^{\varepsilon}\bigg(V_{22}^{\frac{1}{4}}+x^{\frac{53}{54}}\Delta^{\frac{1}{2}}K\bigg)\,.
\end{equation}
The sum (\ref{W22pred}) we will estimate in two ways. If $\Delta<D\le x^{\frac{8}{27}}\Delta$
then from inequality (\ref{Ocz2}) of Lemma \ref{M3S2J} follows
\begin{equation}\label{ocenkaW22HeathBrown}
  V_{22}^{\frac{1}{4}}\ll x^{\varepsilon}\bigg(x^{\frac{205}{216}}\Delta^{\frac{1}{4}}K +\frac{xK}{\Delta^{\frac{1}{16}}}+\frac{xK}{q^{\frac{1}{16}}}
  +x^{\frac{7}{8}}q^{\frac{1}{16}}K^{\frac{15}{16}} \bigg)\,.
\end{equation}
If $D\le\Delta$
we will estimate the sum $V_{22}$ of (\ref{W22pred}) by putting $u=d^3z$. Applying
Lemma \ref{Trsum2} and Lemma \ref{Tau} (iv) we find
\begin{equation}\label{ocenkaW2HeathBrownDelta}
  V_{22}^{\frac{1}{4}}\ll x^{\varepsilon}\bigg(x^{\frac{7}{8}}\Delta^{\frac{3}{4}}K +
\frac{x\Delta ^{\frac{1}{2}}K}{q^{\frac{1}{4}}}
  +x^{\frac{1}{2}}\Delta^{\frac{1}{2}}K^{\frac{3}{4}}q^{\frac{1}{4}}\bigg)\,.
\end{equation}
From (\ref{W2crezW22}), (\ref{ocenkaW22HeathBrown}), (\ref{ocenkaW2HeathBrownDelta}) and (\ref{us1})
follows
\begin{equation}\label{ocenkaW2HeathBrownokonch}
  W_{2}\ll x^{\varepsilon}\bigg(x^{\frac{53}{54}}\Delta^{\frac{1}{2}}K+
\frac{x\Delta ^{\frac{1}{2}}K}{q^{\frac{1}{4}}}+\frac{xK}{q^{\frac{1}{16}}}+\frac{xK}{\Delta^{\frac{1}{16}}}
  +x^{\frac{7}{8}}q^{\frac{1}{16}}K^{\frac{15}{16}}+x^{\frac{1}{2}}\Delta^{\frac{1}{2}}K^{\frac{3}{4}}q^{\frac{1}{4}}
 \bigg)\,.
\end{equation}
It can be seen that for this estimate to be non-trivial to the constraints (\ref{Ogr_q1})
on $q$ we must also impose the constraint
\begin{equation}\label{Ogr_q}
   q<\frac{x^2}{\Delta^2K^3}\,.
\end{equation}

\subsubsection{Evaluation of type I sums.}
$\\$
Let
\begin{equation*}
  D<x^{8/27}\Delta ,\quad M\le x^{1/3},\quad x^{2/3}\le L\,.
\end{equation*}
Reasoning as in the estimation of the sum $W_{13}$ with $H=D$ (see (\ref{W13})) we obtain
\begin{multline}\label{W2malkoDelta}
   W_{1}^2\ll x^{\varepsilon}\bigg(MDK_0\sum\limits_{z\sim MD}\tau (z)
\mathop{\sum}_{\substack{
                      t\ll \frac{LK_0}{D} }}\tau_3(t)
\min\bigg\{\frac{x^2K_0}{z^2t},\,\frac{1}{||\alpha z^2t||}
\bigg\}+xMDK^2\bigg)\\
=x^{\varepsilon}\bigg(V_1+xMDK^2\bigg)\,.
\end{multline}
If $MD>\Delta$ we evaluate sum $V_1$ using Lemma \ref{Mat}:
\begin{multline}\label{ocW1malkoD}
   V_{1}\ll x^{\varepsilon}\bigg(xMDK^2+\frac{x^2K^2}{(MD)^{\frac{1}{2}}}+\frac{x^2K^2}{q^{\frac{1}{2}}}+
xK^{\frac{3}{2}}q^{\frac{1}{2}}\bigg)\\
\ll x^{\varepsilon}\bigg(x^{\frac{44}{27}}\Delta K^2+\frac{x^2K^2}{\Delta^{\frac{1}{2}}}+\frac{x^2K^2}{q^{\frac{1}{2}}}+
xK^{\frac{3}{2}}q^{\frac{1}{2}}
\bigg)
\end{multline}
If $MD\le \Delta$ we put $u=z^2t$ in sum $V_1$ and from Lemma \ref{Trsum2}, Lemma \ref{Tau} (iv) and (\ref{OgrK0})
obtain
\begin{align}\label{ocW1malkoDvtorinachin}
   V_{1}&\ll x^{\varepsilon}MDK_0\sum\limits_{u\ll xMDK_0}\tau_5(u)
\min\bigg\{\frac{x^2K_0}{u},\,\frac{1}{||\alpha u||}
\bigg\}\notag\\
&\ll x^{\varepsilon}\bigg (MDKq+xM^2D^2K^2+\frac{x^2MDK^2}{q}
\bigg)\notag\\
&\ll
x^{\varepsilon}\bigg (\Delta Kq+x\Delta^2K^2+\frac{x^2K^2\Delta}{q}
\bigg)
\end{align}
From (\ref{W2malkoDelta}), (\ref{ocW1malkoD}) and  (\ref{ocW1malkoDvtorinachin}) we receive
\begin{equation}\label{ocW1pyrwa}
   W_{1}\ll x^{\varepsilon}\bigg(\frac{xK}{\Delta^{\frac{1}{4}}}+\frac{xK}{q^{\frac{1}{4}}}+
x^{\frac{1}{2}}K^{\frac{3}{4}}q^{\frac{1}{4}}+\frac{x\Delta^{\frac{1}{2}}K}{q^{\frac{1}{2}}}+\Delta^{\frac{1}{2}}K^{\frac{1}{2}}q^{\frac{1}{2}}
\bigg)\,
\end{equation}
and to the constraints (\ref{Ogr_q1}), (\ref{Ogr_q}) we must also add the constraint
\begin{equation}\label{Ogr_q2}
  q>\Delta K^2\,.
\end{equation}
Using (\ref{ocenkaW2HeathBrownokonch}), (\ref{ocW1pyrwa}) and (\ref{Ogr_q2}) for case $D\le x^{\frac{8}{27}}\Delta $ the estimate
\begin{equation}\label{ocenkaWHeathBrownokonch}
  W\ll x^{\varepsilon}\bigg(
\frac{x\Delta ^{\frac{1}{2}}K}{q^{\frac{1}{4}}}+\frac{xK}{q^{\frac{1}{16}}}+\frac{xK}{\Delta^{\frac{1}{16}}}
  +x^{\frac{7}{8}}q^{\frac{1}{16}}K^{\frac{15}{16}}+x^{\frac{1}{2}}\Delta^{\frac{1}{2}}K^{\frac{3}{4}}q^{\frac{1}{4}}
 \bigg)\,.
\end{equation}
is true.
Form (\ref{ocenkaW3}) and (\ref{ocenkaWHeathBrownokonch})
we get
\begin{equation*}
  W\ll x^{\varepsilon}\bigg(\frac{xK}{\Delta^{\frac{1}{32}}}+\frac{xK}{q^{\frac{1}{32}}}+
  x^{\frac{15}{16}}K^{\frac{31}{32}}q^{\frac{1}{32}} +\frac{x\Delta ^{\frac{1}{2}}K}{q^{\frac{1}{4}}}+
 x^{\frac{1}{2}}\Delta ^{\frac{1}{2}}K^{\frac{3}{4}}q^{\frac{1}{4}} \bigg)\,.
\end{equation*}
It is not difficult to see that for
\begin{equation}\label{Ogr_qOkonch}
  \max\{x^{\omega},\,\Delta K^2\}<q<\min\bigg\{\frac{x^2}{K^{31}},\,\frac{x^2}{\Delta ^2K^{3}}\bigg\}
\end{equation}
a non-trivial estimate
\begin{equation*}\label{ocWokonch}
   W\ll x^{\varepsilon}\bigg(\frac{xK}{\Delta^{\frac{1}{32}}}+
\frac{x\Delta^{\frac{1}{2}}K}{q^{\frac{1}{4}}}
+x^{\frac{15}{16}}K^{\frac{31}{32}}q^{\frac{1}{32}}+
x^{\frac{1}{2}}\Delta ^{\frac{1}{2}}K^{\frac{3}{4}}q^{\frac{1}{4}}
\bigg)\,.
\end{equation*}
is valid.

\subsection{Proof of Lemma \ref{Lemseq}}
$\\$

To prove Lemma \ref{Lemseq} it is enough
in Theorem \ref{Mat2}
choose
\begin{equation}\label{deltaeta}
   x=q,\quad \Delta=K^{\frac{32.34}{33}},\quad K=x^{\theta}\log ^2x,\quad D=\frac{x^{1/2}}{\Delta K^{4}}\,,
\end{equation}
where $\theta<0,005$ is arbitrary small fixed number.

\section{Proof of Theorem \ref{osnTh}}
As in \cite{TT} we take a periodic with period 1 function such that
\begin{equation}\label{hi}
\begin{aligned}
    0<\chi (t) &<1 \quad \mbox { if }\quad  -\delta< t< \delta;\\
    \chi (t) &=0 \quad \mbox { if }\quad\quad\;  \delta \le t\le 1-\delta,
\end{aligned}
\end{equation}
and which has a Fourier series
\begin{equation}\label{hi1}
    \chi (t)= \delta +\sum\limits_{|k|>0}c(k)e(kt),
\end{equation}
with coefficients satisfying
\begin{equation}\label{hi2}
    c(0)=\delta,\quad c(k)\ll \delta \quad \mbox{ for all } k,\quad
    \sum\limits_{|k|>K}|c(k)|\ll    x^{-1}
\end{equation}
and $\delta $ and $K$ satisfying the conditions (\ref{us2}).
The existence of such a function is a consequence of a well known
lemma of Vinogradov (see \cite{Kar}, ch. 1, \S 2).

Next we will use the linear sieve due to Iwaniec - this is Theorem \ref{linSieve} (see \cite{Iwa1}).
As usual, for any sequence $\mathcal{A}$ of integers weighted
by the numbers $f_n$, $n\in \mathcal{A}$, we set
\begin{align*}
  &S(\mathcal{A},\,z)=\sum\limits_{\substack{
                      n\in \mathcal{A}\\
                      (n,\,P(z))=1 }}f_n,\quad\quad P(z)=\prod\limits_{p<z}p,\\
&V(z)=\prod\limits_{p|P(z)}\bigg(1-\frac{\omega (p)}{p}\bigg),\quad
C_0=\prod\limits_{p>2}\bigg(1-\frac{1}{(p-1)^2}\bigg)
\end{align*}
and denote by $\mathcal{A}_d$ be the subsequence of elements $n\in \mathcal{A}$
with $n\equiv 0\pmod d$.
To prove Theorem \ref{osnTh}, it suffice to show that
\begin{equation}\label{Saineq}
  S(\mathcal{A},\,N^{1/3})=\mathop{\sum}_{\substack{
                      p+2\le x\\
                      (p+2,\,P(x^{1/3})=1 }}\chi (\alpha p^2+\beta)>0\,.
\end{equation}
Following the exposition in Shi's article (see \cite{Shi}) we have that
\begin{multline*}
  S\ge  S(\mathcal{A},\,x^{1/12})-\frac{1}{2}\sum\limits_{x^{1/12}\le p_1 < x^{1/3.1}}S(\mathcal{A}_{p_1},\,x^{1/12})
-\frac{1}{2}\sum\limits_{\substack{
                      x^{1/12}\le p_1 < x^{1/3.1}\\
                     x^{1/3.1}\le p_2<(\frac{x}{p_1})^{1/2} }}S(\mathcal{A}_{p_1p_2},\,x^{1/12})
\\
-\sum\limits_{\substack{x^{1/12}\le p_1 < p_2<(\frac{x}{p_1})^{1/2} }}S(\mathcal{A}_{p_1p_2},\,x^{1/12})
+O\big(x^{11/12}\big)\\
=S_1-\frac{1}{2}S_2-\frac{1}{2}S_3-S_4+O\big(x^{11/12}\big)
\end{multline*}
and it is enough to proof that above expression is positive.
Consider a square-free number $d$. If $2|d$, then we write $|\mathcal{A}_d|=|r(\mathcal{A}, d)|\le 1$.
Otherwise we have by the Fourier expansion of $\chi (n)$ that
\begin{align*}
  |\mathcal{A}_d|= &\sum\limits_{\substack{
                      p\le x-2\\
                     p\equiv -2 (d) }}\chi (\alpha p^2+\beta))\\
  = & \sum\limits_{\substack{
                      p\le x-2\\
                     p\equiv -2 (d) }}\bigg(
\delta+\delta\sum\limits_{0<|k|<K}c(k)e(\alpha n^2k)+O(x^{-1})
\bigg)\\
 = & \delta \Bigg(\frac{Li\, x}{\varphi (d)}+R_1(d)+R_2(d)+O\bigg(\frac{1}{d}\bigg)
\Bigg)\,,
\end{align*}
where $A>0$ is arbitrarily large.
Here we laid $c(k):=c(k)e(\beta k)$, $c(k)\ll 1$ and
\begin{equation*}
  R_1(d)=  \sum\limits_{\substack{
                      p\le x-2\\
                     p\equiv -1 (d) }}1- \frac{Li\,x}{\varphi (d)},\quad\quad
  R_2(d)= \sum\limits_{\substack{|k|<K\\k\ne 0}}c(k)\sum\limits_{\substack{
                      p\le x-2\\
                     p\equiv -2 (d) }}e(\alpha p^2k)\,.
\end{equation*}
Applying Bombieri–Vinogradov theorem (see \cite{IwaKow}, Theorem 17.1) we have
\begin{equation*}
  \sum\limits_{d\le D}|R_1(d)|\ll \frac{x}{(\log x)^A}\,.
\end{equation*}
To evaluate the sum 
\begin{equation*}
  \sum\limits_{d\le D}\lambda (d) R_2(d)
\end{equation*}
it is enough to evaluate the sum
\begin{equation*}
  \sum\limits_{d\le D}\lambda (d)\sum\limits_{\substack{|k|<K\\k\ne 0}}c(k)
\sum\limits_{\substack{n\le x-2\\n\equiv -2 (d) }}\Lambda (n)e(\alpha n^2k)\,.
\end{equation*}
On the other hand, Lemma \ref{Lemseq} implies that for a well-separable numbers $d$ of level $D=\dfrac{x^{1/2}}{\Delta K^4}$
and $\lambda (d)\ll \tau (d)$ we get
\begin{equation*}
  \sum\limits_{d\le D}\lambda (d)\sum\limits_{\substack{|k|<K\\k\ne 0}}c(k)
\sum\limits_{\substack{n\le x-2\\n\equiv -2 (d) }}\Lambda (n)e(\alpha n^2k)\ll x^{1-\omega}
\end{equation*}
Here $x=q$, where $a/q$ convergent to $\alpha $ with a large enough denominator.
From here on, the reasoning we go through is the same as in Shi's paper (see \cite{Shi}).
We will only note that to estimate the sum
\begin{equation*}
\sum\limits_{x^{1/12}\le p<x^{1/3.1}}\sum\limits_{d\le D}R_2(pd)
\end{equation*}
with $D=\dfrac{x^{1/2}}{\Delta pK^4}$
first we present it as a $O(\log^4 x)$ number of sums of the type
\begin{equation*}
   R_2(P)=\sum\limits_{d\sim D}\lambda(d)
    \sum\limits_{1\le |k|\sim K}c(k)
\sum\limits_{p\sim P}
\sum\limits_{\substack{ n\sim x \\ n+2\equiv 0\,(d)\\
                                       n+2\equiv 0\, (p_1)}}\Lambda (n)e\big(\alpha n^2k\big )\,,
\end{equation*}
where $x^{\frac{1}{12}}\le P<\dfrac{1}{2}.x^{\frac{1}{3.1}}$.

If $DP\le x^{8/27}\Delta$ we put $t=dp$ and represent the sum $R_2(P)$ in type:
\begin{equation*}
  R_2(P)=\sum\limits_{1\le |k|\sim K}c(k)
\sum\limits_{t\sim DP}g(t,\,d)
\sum\limits_{\substack{ n\sim x \\ n+2\equiv 0\,(t)}}\Lambda (n)e\big(\alpha n^2k\big )\,,
\end{equation*}
where
\begin{equation*}
  g(t,\,d)=\sum\limits_{\substack{ d\sim D \\ d|(t,\,P(z))\\ t/d>x^{1/12}\\t/d - \hbox{\tiny{prime}}}}\lambda(d)\ll \tau (t)
\end{equation*}
and evaluation is in the same way as in \S $\,$\ref{heath}.

If $DP\ge x^{8/27}\Delta$ then, depending on which interval $P$ falls into
and the fact that $d$ is well-separated, we choose $H$ so that $PH$ falls
into one of the intervals $x^{2/5}\ll PH\ll \dfrac{x^{1/2}}{\Delta K^4}$ or
$x^{8/27}\Delta\ll PH\ll x^{2/5}$. So
\begin{multline*}
  R_2(P)=\sum\limits_{1\le |k|\sim K}c(k)
\sum\limits_{s\sim S}\sum\limits_{h\sim H}\lambda(hs)
\sum\limits_{p_1\sim P}
\sum\limits_{\substack{ n\sim x \\ n+2\equiv 0\,(d)\\
                                       n+2\equiv 0\, (p_1)}}\Lambda (n)e\big(\alpha n^2k\big )\\
=\sum\limits_{1\le |k|\sim K}c(k)
\sum\limits_{s\sim S}\sum\limits_{t\sim PH}g(t,\,s)
\sum\limits_{\substack{ n\sim x \\ n+2\equiv 0\,(ts)}}\Lambda (n)e\big(\alpha n^2k\big )
\end{multline*}
where
\begin{equation*}
  g(t,\,s)=\sum\limits_{\substack{ h\sim H \\ h|(t,\,P(z))\\t/h>x^{1/12}\\ t/h-\hbox{\tiny{prime}}}}\lambda(hs)\ll \tau (t)
\end{equation*}
and evaluation is in the same way as in \S $\,$\ref{Vaughan}.
For the sum $S_3$ 
and $S_4$ we use the same calculation as in \cite{Shi}
and choosing
\begin{equation*}
   z=x^{\frac{1}{12}},\quad \Delta=K^{\frac{32.34}{33}},\quad K=x^{\frac{1}{1296}-\eta},\quad \quad D=\frac{x^{1/2}}{\Delta K^{4}}\,,
\end{equation*}
we get that the inequality (\ref{Saineq}) is true and the proof of Theorem
is complete.

\bigskip
\bigskip

\vbox{
\hbox{Faculty of Mathematics and Informatics}
\hbox{Sofia University ``St. Kl. Ohridsky''}
\hbox{5 J.Bourchier, 1164 Sofia, Bulgaria}
\hbox{ }
\hbox{tlt@fmi.uni-sofia.bg}}

\bigskip
\bigskip
\bigskip


\medskip






\end{document}